\documentstyle[12pt]{article}
\long\def\comment#1{}

\begin{document}

\title{Early Record of Divisibility and Primality}
\maketitle

\begin{centering}
\author{Subhash Kak}\\
\end{centering}

\begin{abstract}
We provide textual evidence on divisibility and
primality in the ancient Vedic texts of India. Concern with divisibility becomes
clear from the listing of all the fifteen pairs of divisors of the number 720. The total number of pairs of divisors of 10,800 is also given. The motivation behind finding the  divisors was the theory that the number of divisors of a certain periodic process is related to the count associated with some other periodic process.
For example,  720 (days and nights
of the year) has 15 pairs of divisors, and this was related to
the 15 days of the waxing and waning of the moon.
Numbers that have no divisors appeared to have
been used to symbolize the ``transcendent" that is
beyond periodicity and change.

\end{abstract}

\maketitle

\begin{centering}
{\it Dedicated to my mother on her 80th birthday on 23rd March, 2009}
\end{centering}
\section{Introduction}

The Vedas are early texts from India, conservatively dated to a long period beginning in late third or early second millennium BC,  which are either in verse as hymns, or as prose commentaries on ritual or as aphorisms dealing with various branches of knowledge. Because of their obscure style, they have not been studied carefully for their mathematical knowledge. These texts have many ill-understood mathematical allusions and they
speak of large numbers at many places.

The Vedic literature is to be found in several layers starting with the four Vedas (\d{R}gveda, Yajurveda, S\={a}maveda, and Atharvaveda),  Br\={a}hma\d{n}a texts, \={A}ra\d{n}yakas, Upani\d{s}ads, and S\={u}tras. Some of the materials is arranged in collections called Sa\d{m}hit\={a}s.
The Yajurveda (17.2) gives a sequence of
powers of 10 going to $10^{12}$. In the Vedic book called {\em \'{S}atapatha Br\={a}hma\d{n}a} (\'{S}B) [11], there is a sequence (12.3.2) speaking of different successive divisions of the year that amounts to $10,800 \times 15^6$ parts. Elsewhere the number
of stars is given as $1.08 \times 10^7$. Other numbers
are used symbolically.

A famous verse from the I\'{s}a Upani\d{s}ad speaks of
``fullness'' from which ``fullness'' arises and if ``fullness'' is
subtracted from it ``fullness'' remains, which indicates that the Vedic authors
had the intuition of
the mathematical idea of infinity. Elsewhere, there is explicit  mention
of infinity as being uncountable.

  The texts have stories that have a mathematical basis. For example, the Maitr\={a}ya\d{n}\={\i}
   Sa\d{m}hit\={a} 1.5.8 [1],[2] has the story of Manu with ten wives, who have  one, two, three, four, five, six, seven, eight, nine, and ten sons, respectively.   The one son allied with the nine sons, and the two sons allied with the eight, and so on until the five sons were left by themselves. They asked the father for help, and  he gave them each a samidh, or ``oblation-stick," which they used  to defeat all of the other sons.

 Since the ten sons did not ally with anyone, and the pairing of the others, excepting the five left over, is in
 groups of ten, the counting is in the base 10 system. In this mathematical story, the sticks help make the five stronger than the other 50. Perhaps this happens because each stick has a power of 10, and therefore the 5 now have a
 total power of 55 which vanquishes the 50. This could imply knowledge of the place value system if one conjectures that  each oblation-stick is in the higher place value so that 50+5=55.

In \'{S}B 10.4.2, the total number of syllables in the \d{R}gveda is taken to
be 432,000, equal to the number of muh\={u}rtas in 40 years (each muh\={u}rta = 48 minutes).
The Yajurveda and the S\={a}maveda taken together are taken to be another 432,000 syllables.
Several years ago, I showed that many other numbers
in the Vedic texts also have an astronomical basis [3]-[8].
At that time I did not raise the question whether
the Vedic sages were aware of properties of numbers
beyond their use for numerical calculations.

While doing that research I was aware that
an early index to the \d{R}gveda due to \'{S}aunaka [9]
states that the number of words in the book is
153,826, and this number is twice
76,913, which is a  prime number.
Since the mantras are
in verses of two lines, one might expect that
the number 76,913 is the more basic one and
it was deliberately chosen. Likewise, the
total number of verses in the \d{R}gveda is $10522 = 2 \times 5261$, and 5261 is prime.
But these two examples of primality could just
be accidental,
unless one could show that the Vedic authors actually
knew this concept.

It appears that in the Vedic texts one can
distinguish between numbers that are derived
from observed phenomena and others that are ideal, or have an
abstract basis. It is the latter numbers that
are likely to be prime.
It is due to the assumption
of connections, {\em bandhu-}, between the
astronomical, the physical, and
the elemental, which is central to Vedic thought, that many  numbers
are astronomical and related to
the motions of the sun, moon, and the planets.
Other numbers  show up
where the narrative transcends
astronomy or physical structure and
describes the {\em Puru\d{s}a}, the Cosmic
Man, as in yantras. These are  non-astronomical and
may be prime or have large prime factors.

This article presents  evidence supporting an early
tradition of systematic
examination of properties of numbers.  We begin with the background on mathematical knowledge in the Vedic times. Then we present
examples of the counting of all the divisors of a number and the motivation
for finding numbers that had a specified number of divisors. Although we do not
have evidence from the texts speaking directly of the idea of primality,
we suggest that indirect evidence supports the knowledge of this idea.

\section{The Mathematics of the Vedic Times}

 The Br\={a}hma\d{n}as and the \'{S}ulbas\={u}tras [10] give account of early Vedic mathematics.
Apart from concerns of geometry and astronomy the \'{S}atapatha Br\={a}hma\d{n}a (\'{S}B) [11]
 deals with the question of all the
divisors of a number.
It also provides approximations to $\pi$ [12]. For general surveys of Indian mathematics, see [13]-[16].

 The \'{S}ulbas\={u}tras give geometric solutions of linear
 equations in a single unknown. They also deal with quadratic equations of the forms $ax^2 = c$ and $a x^2 + bx = c$.
Baudh\={a}yana's \'{S}ulbas\={u}tra gives a remarkable approximation to $\sqrt 2$:

\vspace{0.2in}

\[\sqrt 2 = 1 + \frac{1}{3} + \frac{1}{3 \times 4} - \frac{1}{3 \times 4 \times 34} = \frac{577}{408} \]

\vspace{0.2in}

This is accurate to five decimal places. It is intriguing that Baudh\={a}yana felt the need to add the last term
in the expansion because without that the approximation is still valid to three decimal places and excellent for geometric constructions.

The motivation for the mathematics of the \'{S}ulbas\={u}tras is the solution to  problems of
altar construction. Some of these constructions are squaring a circle (derived from the equivalence of
the circular earth altar and the square sky altar) and construction of a geometric design of a larger size
by increasing the dimensions. Since in the actual construction of such altars the accuracy of the above expansion
of $\sqrt 2$ would not have been noticed, it is clear that Baudh\={a}yana was interested in mathematical problems and
properties of numbers.

The \'{S}ulbas\={u}tras belong to the Ved\={a}\.{n}gas, or supplementary texts of the Vedas. Although they are part of the Kalpa S\={u}tras, which deal with ritual, their importance stems from the constructions they provide for building geometric altars. Their contents, written in the condensed s\={u}tra style, include geometrical propositions and problems related to rectilinear figures and their combinations and transformations, squaring the circle, as well as arithmetical and algebraic solutions to these problems. The root {\em \'{s}ulb} means measurement, and the word ``\'{s}ulba'' means a cord, rope, or string.

The extant \'{S}ulbas\={u}tras belong to the schools of the Yajurveda. The most important \'{S}ulba texts are the ones by Baudh\={a}yana, \={A}pastamba, K\={a}ty\={a}yana, and M\={a}nava. They have been generally assigned to the period 800 to 500 B.C., although they are likely to be older. Baudh\={a}yana's text is the oldest, and he is believed to have belonged to the Andhra region.
Baudh\={a}yana begins with units of linear measurement and then presents the geometry of rectilinear figures, triangles, and circles, and their transformations from one type to another using differences and combinations of areas. An approximation to the square root of 2 and to $\pi$ are next given.

Then follow constructions for various kinds of geometric altars in the shapes of the falcon (both rectilinear and with curved wings and extended tail), kite, isosceles triangle, rhombus, chariot wheel with and without spokes, square and circular trough, and tortoise.

In the methods of constructing squares and rectangles, several examples of Pythagorean triples are provided. It is clear from the constructions that both the algebraic and the geometric aspects of the so-called Pythagorean theorem were known. This knowledge precedes its later discovery in Greece. The other theorems in the \'{S}ulba include:

\begin{itemize}
	\item The diagonals of a rectangle bisect each other.

	\item The diagonals of a rhombus bisect each other at right angles.

	\item The area of a square formed by joining the middle points of the sides of a square is half of the area of the original one.
	
	\item A quadrilateral formed by the lines joining the 	
middle points of the sides of a rectangle is a rhombus
	whose area is half of that of the rectangle.

	\item A parallelogram and rectangle on the same base and
	within the same parallels have the same area.

	\item If the sum of the squares of two sides of a triangle
	is equal to the square of the third side, then the triangle
	is right-angled.
\end{itemize}

A variety of constructions are listed. Some of the geometric constructions in these texts are based on algebraic solutions of simultaneous equations, both linear and quadratic. It appears that geometric techniques were often used to solve algebraic problems.

The \'{S}ulbas are familiar with fractions. Algebraic equations are implicit in many of their rules and operations. For example, the quadratic equation and the indeterminate equation of the first degree are a basis of the solutions presented in the constructions.

The \'{S}ulba geometry was used to represent astronomical facts. The altars that were built according to the \'{S}ulba rules demonstrated knowledge of the lunar and the solar years.

To consider the larger background of Indian science requires an understanding of
its cosmology [17]-[21].
Seidenberg [22],[23] examined the
geometry and mathematics
of the   \'{S}atapatha Br\={a}hma\d{n}a
and he argued that the philosophy that equivalent altars were to have
equal areas led to the posing of basic problems of geometry
leading to results such as the Pythagoras theorem.
A conservative chronology has placed the
final form of this book to
1000-800 B.C.; new archaeological discoveries suggest that
this book may be a thousand years older.
There is also a specific astronomical reference in the book, namely that
the Pleiades do not swerve from the east, that was true only in late
third millennium B.C. Corroborating evidence comes from the text, since it
claims to belong to a period when the Sarasvati river had recently dried up.
New research shows that this drying up took place around 2000 B.C.

The Vedic Indians were also interested  in metres and music and related
mathematical problems [24], [25]. The earliest codified Vedic
astronomy is given in Lagadha's Ved\={a}\.{n}ga Jyoti\d{s}a [26].

\section{Number Divisors}
The number 360, the days in the year, forms a starting point in the
design of the altars of \'{S}B that were used to represent various
astronomical facts about the year.

\'{S}B 10.4.2.1-18 points out that 720, the nights and days in
the nominal year, has exactly 15 factors that are smaller than
the companion. The text is speaking the year as Praj\={a}pati (the Lord of the people)
represented by 720 bricks who is considering the various divisors of the number [11]:

\begin{quotation}
He divided his body into two: there were three hundred and sixty bricks in the one, and as many
in the other; he did not succeed [\emph{in finding all the divisors}].

He made himself three bodies, and in each there were three times eighty (240) bricks; he did
not succeed.

....

He made himself six bodies of a hundred and twenty brick each; he did succeed.
He did not divide himself sevenfold.

He made himself eight bodies of ninety bricks each; he did not succeed.

....

He made himself twenty bodies of thirty-six brick each; he did not succeed. He did not
divide himself either twenty-one fold, or twenty-two fold, or twenty-three fold.

He made himself twenty-four bodies of thiry bricks each. There he stopped, at the fifteenth,
and because he stopped at the fifteenth arrangement [\emph{having found all the divisors}], there
are fifteen forms of the waxing, and fifteen forms of the waning (moon).
\end{quotation}

These divisors are thus
stated to be: 1, 2, 3, 4, 5, 6, 8, 9, 10, 12, 15, 16,
18, 20, 24.
The numbers in this
sequence which do not divide 720 (that is the numbers 7, 11, 13, 14,
17, 19, 21, 22, 23) are also explicitly mentioned.

The next passage claims that 24, the largest of these numbers,
represents the number of half months in the year, which is
a reiteration of $ 15 \times 24 = 360 $.

This implies the thesis that the moon's circuit
around the earth is 30 days because 720, the number of days and nights
in the year, has exactly 30 different divisors.
In other words, numbers are used to relate a characteristic of the
motion of the sun to another related to the moon's orbit
around the earth. An astronomical fact is ``explained" based
on an abstract numerical property.

It appears that the authors did not rather consider the number 360, which has
12 divisors, which could have also been taken to correspond to the 12 months
of the year, yielding, in turn, the 30-day duration of the month, because the
largest of these divisors which is less than its companion, 18, has no direct
correspondence with the basic facts of the year.

The concern with the number of divisors implies finding out
what numbers do not have divisors, or are prime. We know, that for
number $n = {p_1}^{a_1} \times {p_2}^{a_2} $, the number of divisors $d(n)$ is given by:

\[d(n) = (a_1 +1 ) (a_2 +1) \]

The pairs of these divisors, in the manner of counting by the Vedic authors, is $\delta (n) = d(n)/2$.

The Vedic authors were also interested in the largest divisor, whose companion is smaller than itself.
If this divisor is called $\nu (n)$, we have $\nu (720) = 24$.

\subsection{Divisors of sixty and the 6-day week}
One would expect that paralleling a justification of the thirty day month
based on the number of factors of 720,
an argument would have been used to define divisions of the month.
The argument would look at the 60 days and nights of the month and
determine the number of divisors that,
paralleling the procedure in \'{S}B, are less than the companion.
These divisors are 1, 2, 3, 4, 5, and 6.
This suggests a week of six days.
We do get references to the six-day week called \d{s}a\d{d}aha in the
Vedic texts.
Five \d{s}a\d{d}ahas made a month.
This defined a symmetry with the year of five seasons and the yuga of
five years.

Although it has often been assumed that the seven-day week was a later
innovation, it is quite possible that it was the older tradition and that
the six-day week got mentioned in the texts because of the ``theory''
behind its derivation.
The seven-day week is a more natural,
system  because it divides the lunar
month into four equal parts.

The number 7 plays an important number in Vedic cosmology since it appears
in conjunction with the name of the entire country, Sapta Sindhu, with the
additional ideas of seven rivers, seven continents, seven islands, seven
mountains, seven rishis (the Pleiades), seven musical notes, and seven worlds.

\subsection{The division of the year into 10,800 muh\={u}rtas}
\'{S}B 10.4.2.36 speaks of the division of the year into 10,800 muh\={u}rtas. It is further
stated that the divisors of this number go into 30 arrangements, or 30 pairs of
divisors.   Therefore, $\delta (10,800) = 30$, or  $d(10,800) = 60$. This is indeed correct:

\[d(10,800) = d(2^4 \times 3^3 \times 5^2 ) = (4+1)(3+1)(2+1) = 60 \]

The discovery of the 30 pairs of divisors of 10,800 suggests that this came by searching divisors of several suitable astronomical numbers for one  that had precisely
30 pairs of divisors. This indicates that
the Vedic authors knew that number of divisors varied and that some numbers had no divisors.

The text relates the 30 pairs of divisors of 10,800 to the 30 nights of the month.





\section{Astronomical and Non-astronomical  Numbers}
Astronomical numbers in the Vedic texts are related to
the 360 {\em tithis} of the lunar year,
the 365 or 366 days of the solar year (or 371 or 372
{\em tithis}), the 27 or 28 nak\d{s}atras,
the 29 days of the month,
numbers related to the divisions of the year,
the planet orbits, synchronization periods
for the luni-solar motions or planet motions, and so on.
Another astronomical number is 108, the distance
to the sun or the moon in terms of multiples of their
respective diameters [3], and it figures in the earliest
temple design [27]. This number  is also related to
the nak\d{s}atra year of 324 days (12 ``months'' of 27 lunar days)

Astronomical numbers are generally highly composite.
The best examples are the cosmic cycle numbers,
such as the longest cycle of 311,040,000 million
years.

Other numbers, which  appear non-astronomical,
may actually have an astronomical basis.
For example, the 33 gods are likely related
to the count of 27 nak\d{s}atras, five planets,
and the moon.
Or consider \={A}yurvedic
physiology where the number of joints in
the body ({\it marma}) is taken to be 107, a prime
number.
Its primality may not be the
reason for its choice, because
the 108 disks of the sun from the earth
may have been taken to be mirrored in 108
{\em links} from the feet to the crown,
and these 108 links will then have 107 joints.

Some of the non-astronomical origin numbers may not  be actual
counts, but rather ideal counts and, therefore, the
choice of the number as a prime becomes
significant.

A yantra is a representation of the universe. The
earliest yantras, in the shape of buildings,  have been excavated [28]
in North Afghanistan and they date to 2000 BC.
Some people see in Atharvaveda 10.2.31-32,
which speaks of a eight-wheeled, nine-doored
impregnable stronghold of the gods,
the earliest textual reference to a yantra.
The famous \'{S}ri Yantra, created out of
9 juxtaposed triangles, leads to 43 smaller triangles [29]
Whether the primality of 43 is significant here, we cannot tell.

Another description, that is almost certainly that of
a yantra, is  found in the \'{S}vet\={a}\'{s}vatara
Upani\d{s}ad 1.4: ``The Wheel of Brahman has
one felly, a triple tire, sixteen end-parts, fifty spokes
with twenty counter-spokes, and six sets of eight; whose
one rope is manifold; (which moves on three different roads;
and whose illusion arises from two causes).''
The definition of this yantra is complete with the
description of the rope, the remainder describes two
characteristics of it.
The one rope ({\em p\={a}\'{s}a}) is the
central focus, the {\em bindu} which holds it together;
the six sets of eight appear to be the 6 central triangles,
three of which are pointed upward and three are pointed
downward.
The total count here is $1+1+16+50+20+48+1 = 137$,
a prime number.

\section{A Prime Number from Ideal Physiology}
In Vedic thought, the human body is taken to mirror the
universe, and it is supposed to have cyclic processes that are synchronized
to that of the sun, the moon and other luminaries.
The {\em n\={a}\d{d}\={\i}s}, nerves, are taken to be the
web within the body where impressions are stored.
It is also believed that
by concentration one can transcend these
impressions and achieve union with Brahman.

An early count
of the {\em n\={a}\d{d}\={\i}s}
is given in
the B\d{r}had\={a}ra\d{n}yaka Upani\d{s}ad
2.1.19. It speaks of 72,000 {\em n\={a}\d{d}\={\i}s} called
{\em hit\={a}\d{h}} that spread from
the heart. The Ch\={a}ndogya U. 8.6.6.
speaks of 101 {\em n\={a}\d{d}\={\i}s} from the heart.

The Pra\'{s}na U. 3.6 has the complete
description where it is stated that
there are 101 chief {\em n\={a}\d{d}\={\i}s}, each with
100 branch {\em n\={a}\d{d}\={\i}s}, each of which in
turn has 72,000 tributary {\em n\={a}\d{d}\={\i}s}.

This makes a total of
$101 + 101 \times 100 + 101 \times 100 \times 72,000 = 727,210,101$.
This is equal to $101 \times 7,200,101$, each of which
factor is a prime number.

The prime number 7,200,101
is the number of {\em n\={a}\d{d}\={\i}s} emerging from each
of the 101 chief {\em n\={a}\d{d}\={\i}s} of the body.

\section{Conclusions}

Vedic texts are written in a style that presents
several difficulties to the reader uninitiated in its
cosmology. The reader needs to understand that narratives
on ritual have a subtext that deals with astronomy and/or
physiology because of the assumed equivalences between
the microworld and the macroworld. These
texts  are associated with a numerical mysticism in which
the many (numbers) emerge from the one (number) which is why properties
of numbers were of interest to the Vedic authors.

We have presented evidence that the Vedic
authors connected the number of divisors to certain
periodic processes. The texts provide examples of
systematic calculation of the divisors, suggesting that
they were aware that other numbers do not have divisors
excepting one and the number itself. The matter of counting of
divisors is given prominence in the \'{S}atapatha Br\={a}hma\d{n}a.

Multiplication tables as well as Pythagorean triples have been
obtained from tablets belonging to the Old Babylonian period
of about 1700 BC [30].
But these lists do not present a systematic derivation of any
number-theoretic property. The examination of the evidence from
the Vedas is, therefore, important in  the history of ideas.

The primality of the
 count
of the {\em n\={a}\d{d}\={\i}s},   the
word and verse counts in the \d{R}gveda, and the yantra counts could be
coincidental. But the fact that the Vedic authors were counting divisors
suggests that these numbers were deliberately chosen.

Summarizing, we show that the number of
divisors of numbers such as 720 and 10,800 was systematically
calculated. The concern with  a complete list of divisors of a number suggests
knowledge of primality although it cannot be the proof of that knowledge.

\textbf{Acknowledgement.} I thank A. Raghuram for comments on an earlier version of the paper.

\section*{References}

\begin{enumerate}

\item S.D. Satvalekar (ed.), {\em Maitrayani Samhita.} Government Press, Bombay, 1941.

\item M.D. Pandit, {\em Mathematics as known to the Vedic Samhitas.}  Sri Satguru Publications, Delhi, 1993.

\item S. Kak. The astronomy of the Vedic altars and the \d{R}gveda.
{\em Mankind Quarterly}, 33, 43-55, 1992.

\item S. Kak. The astronomy of the Vedic altars.
{\em Vistas in Astronomy}, 36, 117-140, 1993.


\item S. Kak,  The structure of the \d{R}gveda.
{\em Indian Journal of History of Science,}
28, 71-79, 1993.

\item S. Kak, From Vedic science to Vedanta. {\em Adyar Library Bulletin,}
59, 1-36, 1995.


\item S. Kak, Archaeoastronomy and literature. {\it Current Science,} 73, 624-627, 1997.

\item S. Kak, Time, space and structure in ancient India. Presented at the Conference on Sindhu-Sarasvati Valley Civilization: A Reappraisal, Loyola Marymount University, Los Angeles, February 21 \& 22, 2009; arXiv:0903.3252

\item
A.A. Macdonell, {\it K\={a}ty\={a}yana's Sarv\={a}nukrama\d{n}\={\i}
of the Rigveda.} Clarendon Press, Oxford, 1886.

\item S.N. Sen and A.K. Bag, {\it The \'{S}ulbasutras}. Indian National Science Academy, New Delhi, 1983.

\item J. Eggeling (tr.), {\it The Satapatha Brahmana.} Charles Scribner's Sons, New York, 1900;
http://www.sacred-texts.com/hin/sbr/index.htm

\item S. Kak, Three old Indian values of pi. {\it Indian Journal of History of Science}, 32, 307-314, 1997.

\item B. Datta and A.N. Singh, {\it History of Hindu Mathematics}.  Asia Publishing House, Bombay, 1962.

\item C.N. Srinivasiengar, {\it The History of Ancient Indian Mathematics}.  The World Press, Calcutta, 1967.

\item G.G. Joseph, {\it The Crest of the Peacock, non-European roots of Mathematics}. Princeton University Press, Princeton, 2000.

\item I.G. Pearce, {\it Indian Mathematics: Redressing the balance.} 2002. http://www-groups.dcs.st-and.ac.uk/~history/Projects/Pearce/index.html

\item S. Kak, The Indus tradition and the Indo-Aryans. {\it Mankind Quarterly},  32,
195-213, 1992..

\item S. Kak, Greek and Indian cosmology: review of early history. In {\it The
Golden Chain}. G.C. Pande (ed.). CSC, New Delhi, 2005;
arXiv: physics/0303001.

\item S. Kak, Indian physics: outline of early history. arXiv: physics/0310001.

\item S. Kak, {\it The Gods Within.} Munshiram Manoharlal, New Delhi, 2002.

\item S. Kak, {\it The Architecture of Knowledge.} CSC, Delhi, 2004.

\item A. Seidenberg. The ritual origin of geometry.
{\em Archive for History of
Exact Sciences}, 1, 488-527, 1962.

\item A. Seidenberg. The origin of mathematics. {\em Archive for History of
Exact Sciences}, 18, 301-342, 1978.

\item S. Kak, The golden mean and the physics of aesthetics. {\em Foarm Magazine}, 5, 73-81, 2006; arXiv:physics/0411195


\item K.D. Dvivedi and S.L. Singh, {\em The Prosody of Pingala.} Vishwavidyalaya Prakashan, Varanasi, 2008.

\item S. Kak,  The astronomy of the age of geometric altars.
{\em Quarterly Journal of the Royal Astronomical Society,} 36, 385-396, 1995.

\item S. Kak, The axis and the perimeter of the temple. {\em Kannada Vrinda Seminar Sangama} 2005 held at Loyola Marymount University in Los Angeles, November 19, 2005; arXiv:0902.4850

\item V.I. Sarianidi, {\em Drevnie zemeldel'tsy Afganistana.}
Moskva, 1977.
\item Kulaichev, A.P., ``\'{S}riyantra and its mathematical properties,''
{\em Indian Journal of History of Science}, 19,  279-292, 1983.

\item O. Neugebauer. {\em The Exact Sciences In Antiquity}.
Dover Publications, New York, 1969.

\end{enumerate}

\end{document}